\documentclass[conference]{IEEEtran}
\IEEEoverridecommandlockouts

\usepackage{cite}
\usepackage{graphicx}
\usepackage[cmex10]{amsmath}
\usepackage{amssymb}
\usepackage[boxed]{algorithm2e}
\newcommand{\tr}{ {^\top} }

\newcommand{\set}[1]{ {\mathcal{#1}} }

\newcommand{\expect}[1]{ {\mathbb{E}\left[ {#1} \right]} }
\newcommand{\prob}[1]{ {\mathbb{P}\left[ {#1} \right]} }

\newcommand{\abs}[1]{ {\left| {#1} \right|} }

\newcommand{\norm}[1]{ {\left\lVert {#1} \right\rVert} }
\newcommand{\chull}[1]{ {\overline{ {#1} }} }
\newcommand{\tavg}[1]{ {\overline{ {#1} }} }
\newcommand{\indicator}[1]{ {\mathbb{I}\left[ {#1} \right]} }
\newcommand{\RealSet}{ {\mathbb{R}} }

\newcommand{\defeq}{ {\triangleq} }

\newcommand{\arginf}{\operatornamewithlimits{arginf}}

\newcommand{\maximize}{ {\text{Maximize}} }
\newcommand{\minimize}{ {\text{Minimize}} }
\newcommand{\subjectto}{ {\text{Subject to}} }

\newcommand{\prts}[1]{ {\left[ {#1} \right]} }
\newcommand{\prtr}[1]{ {\left( {#1} \right)} }
\newcommand{\prtc}[1]{ {\left\{ {#1} \right\}} }

\newcommand{\optvar}{ {(\text{opt})} }

\newcommand{\random}{ {\omega} }
\newcommand{\setrandom}{ {\Omega} }
\newcommand{\ssp}{ {\pi} }
\newcommand{\drift}{ {\Delta} }
\newcommand{\lagran}{ {\Gamma} }
\newcommand{\dvar}{ {\lambda} }

\newcommand{\poly}{ {\text{P}} }
\newcommand{\gen}{ {\text{G}} }

\newtheorem{theorem}{Theorem}
\newtheorem{lemma}{Lemma}
\newtheorem{assumption}{Assumption}

\begin{document}

\title{Time-Average Stochastic Optimization with Non-convex Decision Set and its Convergence}

\author{\IEEEauthorblockN{Sucha Supittayapornpong,~~Michael J. Neely}
  \IEEEauthorblockA{Department of Electrical Engineering\\
    University of Southern California\\
    Los Angeles, California\\
    Email: supittay@usc.edu,~~mjneely@usc.edu}
%  \thanks{S. Supittayapornpong and M. J. Neely are with Electrical Engineering Department, University of Southern California, 3740 McClintock Ave., Los Angeles, CA 90089-2565, {\tt\small supittay@usc.edu, mjneely@usc.edu}}%
  %\thanks{This material was presented in part at the IEEE 
  %  International Conference on Communications, Ottawa, 
  %  Canada, June 2012 \cite{Sucha:TwoHop}.

  %  This material is supported in part by one or more of
  %  the following: the NSF Career grant CCF-0747525, the
  %  Network Science Collaborative Technology Alliance sponsored
  %  by the U.S. Army Research Laboratory W911NF-09-2-0053.}
}


\maketitle

\begin{abstract}
This paper considers \emph{time-average stochastic optimization}, where a time average decision vector, an average of decision vectors chosen in every time step from a time-varying (possibly non-convex) set, minimizes a convex objective function and satisfies convex constraints.  This formulation has applications in networking and operations research.  In general, time-average stochastic optimization can be solved by a Lyapunov optimization technique.  This paper shows that the technique exhibits a transient phase and a steady state phase.  When the problem has a unique vector of Lagrange multipliers, the convergence time can be improved.  By starting the time average in the steady state the convergence times become $O(1/\epsilon)$ under a locally-polyhedral assumption and $O(1/\epsilon^{1.5})$ under a locally-non-polyhedral assumption, where $\epsilon$ denotes the proximity to the optimal objective cost.  Simulations suggest that the results may hold more generally without the unique Lagrange multiplier assumption.
\end{abstract}

\section{Introduction}
Stochastic network optimization can be used to design dynamic algorithms that optimally control communication networks \cite{Neely:SNO}.  The framework has several unique properties which do not exist in a traditional convex optimization setting.  In particular, the framework allows for a time-varying and possibly non-convex decision set.  For example, it can treat a packet switch that makes binary $(0/1)$ scheduling decisions, or a wireless system with randomly varying channels and decision sets.

This paper considers \emph{time-average stochastic optimization}, which is useful for example problems of network utility maximization \cite{Neely:TON2008, NOW, Stolyar:Greedy, Eryilmaz:QLB}, energy minimization \cite{Neely:IT2006, Lin:Energy}, and quality of information maximization \cite{Sucha:QoI}.

Time $t \in \prtc{0, 1, 2, \dotsc}$ is slotted.  Define $\setrandom$ to be a finite or countable infinite sample space of random states.  Let $\random(t) \in \setrandom$ denote a random state at time $t$.  Random state $\random(t)$ is assumed to be independent and identically distributed (i.i.d.) across time slots.  The steady state probability of $\random \in \setrandom$ is denoted by $\ssp_\random$.  Let $I$ and $J$ be any positive integers.  Each slot $t$, decision vector $x(t) = (x_1(t), \dotsc, x_I(t))$ is chosen from a decision set $\set{X}_{\random(t)}$.  Define $\tavg{x}$ as the time average vector:
\begin{equation*}
  \tavg{x}(T) \defeq \frac{1}{T} \sum_{t=0}^{T-1} \expect{ x(t) }.
\end{equation*}
The goal is to make decisions over time to solve:
\begin{align}
  \label{problem:sto}
  \minimize \quad & \limsup_{T \rightarrow \infty} f( \tavg{x}(T) ) \\
  \subjectto \quad & \limsup_{T \rightarrow \infty} g_j( \tavg{x}(T) ) \leq 0 && j \in \prtc{1, \dotsc, J} \notag\\
  & x(t) \in \set{X}_{\random(t)} && t \in \prtc{0, 1, 2, \dotsc}. \notag
\end{align}
Here it is assumed that $\set{X}_\random$ is a closed and bounded subset of $\RealSet^I$ for each $\random \in \setrandom$.  Let $\set{C}$ be a closed and bounded set that contains $\cup_{\random \in \setrandom}\set{X}_{\random}$, and $f$ and $g_j$ are convex functions from $\chull{\set{C}}$ to $\RealSet$, where $\chull{\set{A}}$ denotes a convex hull of set $\set{A}$.  Results in \cite{Neely:SNO} imply that the optimal point can be achieved with an ergodic policy for which the limiting time average expectation exists.

Solving formulation \eqref{problem:sto} using the stochastic network optimization framework does not require any statistical knowledge of the random states.  However, if the steady state probabilities are known, the optimal objective cost of formulation \eqref{problem:sto} is identical the optimal cost of the following problem:
\begin{align}
  \label{problem:det}
  \minimize \quad & f(x) \\
  \subjectto \quad & g_j(x) \leq 0 && j \in \prtc{1, \dotsc, J} \notag\\
  & x \in \tavg{\set{X}}, \notag
\end{align}
where $\tavg{\set{X}} \defeq \sum_{\random \in \setrandom} \ssp_\random \tavg{\set{X}}_\random$.  Note that $\alpha \set{A} + \beta \set{B} = \prtc{ \alpha a + \beta b: a \in \set{A}, b \in \set{B}}$ for any $\alpha, \beta \in \RealSet$ and any sets $\set{A}$ and $\set{B}$.

Formulation \eqref{problem:det} is convex; however, its optimal solution may not be in any of the sets $\set{X}_\random$.  In fact, determining whether $x$ is a member of $\tavg{\set{X}}$ may already be a difficult task.  This illustrates that traditional and state-of-the-art techniques for solving convex optimization cannot be applied directly to solve problem \eqref{problem:sto}.  Nevertheless, their convergence times are compelling to be mentioned for a purpose of comparison.

The convergence time of an algorithm is usually measured as a function of an $\epsilon$-approximation to the optimal solution.  For a convex optimization problem, several techniques utilizing a time average solution \cite{Nesterov:dual_averaging,Nedic:approximate_primal,Neely:DistComp} have $O(1/\epsilon^2)$ convergence time.  For unconstrained optimization without restricted strong convexity property, the optimal first-order method \cite{Nesterov:O,Tseng:Fast_gradient} has $O(1/\sqrt{\epsilon})$ convergence time, while the gradient and subgradient methods have respectively $O(1/\epsilon)$ and $O(1/\epsilon^2)$ convergence time \cite{Boyd:CO}.  For constrained optimization, two algorithms developed in \cite{Ermin:AsynADMM,Beck:FastNUM} have $O(1/\epsilon)$ convergence time; however, the results rely on special structures of their formulation.  All of these results are for convex optimization problems, which is not formulation \eqref{problem:sto}.

This paper considers a drift-plus-penalty algorithm, developed in \cite{Neely:SNO}, that solves formulation \eqref{problem:sto}.  The algorithm is shown to have $O(1/\epsilon^2)$ convergence time in \cite{Neely:Convergence}.  Note that a deterministic version of formulation \eqref{problem:sto} and its corresponding algorithm are studied in \cite{Sucha:TavgDet}.

Inspired by the analysis in \cite{Longbo:Delay_reduction}, the drift-plus-penalty algorithm is shown to have a \emph{transient phase} and a \emph{steady state phase}.  These phases can be analyzed in two cases that depend on the structure of a dual function.  The first case is when a dual function satisfies a \emph{locally-polyhedral} assumption and the expected transient time is $O(1/\epsilon)$.  The second case is when the dual function satisfies a \emph{locally-non-polyhedral} assumption and the expected transient time is $O(1/\epsilon^{1.5})$.  Then, under a uniqueness assumption on Lagrange multipliers, if the time average starts in the steady state, a solution converges in $O(1/\epsilon)$ and $O(1/\epsilon^{1.5})$ time slots under the locally-polyhedral and locally-non-polyhedral assumptions respectively.  Simulations suggest these results may hold more generally without the uniqueness assumption.

The paper is organized as follows.  Section \ref{sec:dpp} constructs an algorithm solving problem \eqref{problem:sto}.  The behavior and properties of the algorithm are analyzed in Section \ref{sec:behavior}.  Section \ref{sec:poly} analyzes the transient phase and the steady state phase under the locally-polyhedral assumption.  Results under the locally-non-polyhedral assumption are provided in Section \ref{sec:nonpoly}.  Simulations are performed in Section \ref{sec:simulation}.

\section{Time-Average Stochastic Optimization}
\label{sec:dpp}
A solution to problem \eqref{problem:sto} can be obtained through an auxiliary problem, which is formulated in such a way that its optimal solution is also an optimal solution to the time-average problem.  To formulate this auxiliary problem, an additional set and mild assumptions are defined.  First of all, it can be shown that $\chull{\set{X}}$ is compact.

%\begin{assumption}
%  \label{ass:Xbound}
%  Set $\tavg{\set{X}}$ is closed and bounded.
%\end{assumption}

\begin{assumption}
  \label{ass:slater}
There exists a vector $\hat{x}$ in the interior of $\tavg{\set{X}}$ that satisfies $g_j(\hat{x}) < 0$ for all $j \in \prtc{1, \dotsc, J}$.
\end{assumption}
In convex optimization, Assumption \ref{ass:slater} is a \emph{Slater condition}, which is a sufficient condition for strong duality \cite{Bertsekas:Convex}.

Define the \emph{extended set} $\set{Y}$ that is a closed, bounded, and convex subset of $\RealSet^I$ and contains $\tavg{\set{X}}$.  Set $\set{Y}$ can be $\tavg{\set{X}}$, but it can be defined as a hyper-rectangle set to simplify a later algorithm.  Define $\norm{ \cdot }$ as the Euclidean norm.

\begin{assumption}
  \label{ass:lipschitz}
Functions $f$ and $g_j$ for $j \in \prtc{1, \dotsc, J}$ are convex and \emph{Lipschitz continuous} on the extended set $\set{Y}$, so there are constants $M_f > 0$ and $M_{gj} > 0$ for $j \in \prtc{1, \dotsc, J}$ that for any $x, y \in \set{Y}$:
\begin{align}
  \abs{ f(x) - f(y) } & \leq M_f \norm{ x - y } \label{eq:f_lip} \\
  \abs{ g_j(x) - g_j(y) } & \leq M_{gj} \norm{ x - y }. \label{eq:gj_lip}
\end{align}
\end{assumption}

\subsection{Auxiliary formulation}
For function $a(x(t))$ of vector $x(t)$, define an average of function values as
\begin{equation*}
  \tavg{a(x)} \defeq \lim_{T \rightarrow \infty} \frac{1}{T} \sum_{t = 0}^{T-1} \expect{ a(x(t)) }.
\end{equation*}

Recall that problem \eqref{problem:sto} can be achieved with an ergodic policy for which the limiting time average expectation exists.  The time average stochastic optimization \eqref{problem:sto} is solved by considering an \emph{auxiliary formulation}, which is formulated in terms of well defined limiting expectations ``for simplicity.''
\begin{align}
  \label{problem:aux}
  \minimize \quad & \tavg{ f(y) } \\
  \subjectto \quad & \tavg{ g_j(y) } \leq 0 && j \in \prtc{1, \dotsc, J} \notag\\
  & \lim_{T \rightarrow \infty} \tavg{x}_i(T)  =  \lim_{T \rightarrow \infty} \tavg{y}_i(T) && i \in \prtc{1, \dotsc, I} \notag\\
  & (x(t), y(t)) \in \set{X}_{\random(t)} \times \set{Y} && t \in \prtc{0, 1, 2, \dotsc}. \notag
\end{align}
This formulation introduces the auxiliary vector $y(t)$.  The second constraint ties $\lim_{T \rightarrow \infty} \tavg{x}(T)$ and $\lim_{T \rightarrow \infty} \tavg{y}(T)$ together, so the original objective function and constraints of problem \eqref{problem:sto} are preserved in problem \eqref{problem:aux}.  Let $f^\optvar$ be the optimal objective cost of problem \eqref{problem:sto}.

\begin{theorem}
The time-average stochastic problem \eqref{problem:sto} and the auxiliary problem \eqref{problem:aux} have the same optimal cost, $f^\optvar$.

\begin{IEEEproof}
Let $\hat{f}^\optvar$ be the optimal objective cost of the auxiliary problem \eqref{problem:aux}.  We show that $\hat{f}^\optvar = f^\optvar$.

Let $\prtc{ x^\ast(t) }_{t=0}^\infty$ be an optimal solution, generated by an ergodic policy, to problem \eqref{problem:sto} such that:
\begin{align*}
  \lim_{T \rightarrow \infty} f( \tavg{ x^\ast } (T) ) & = f^\optvar \\
  \lim_{T \rightarrow \infty} g_j( \tavg{ x^\ast } (T) ) & \leq 0 && j \in \prtc{1, \dotsc, J} \\
  x^\ast(t) & \in \set{X}_{\random(t)} && t \in \prtc{0, 1, 2, \dotsc }.
\end{align*}
Consider a solution $\prtc{ x(t), y(t) }_{t=0}^\infty$ to problem \eqref{problem:aux} as follows:
\begin{equation*}
  x(t) = x^\ast(t) \quad\quad y(t) = \lim_{T \rightarrow \infty} \tavg{ x^\ast }(T) \quad\quad t \in \prtc{0, 1, 2, \dotsc}.
\end{equation*}
It is easy to see that this solution satisfies the last two constraints of problem \eqref{problem:aux}.  For the first constraint, it follows from Lipschitz continuous that for $j \in \prtc{1, \dotsc, J}$
\begin{equation*}
  \tavg{ g_j(y) } = \tavg{ g_j( \lim_{T \rightarrow \infty} \tavg{x^\ast}(T)) } = g_j( \lim_{T \rightarrow \infty} \tavg{x^\ast}(T)) \leq 0.
\end{equation*}
Therefore, this solution is feasible, and the objective cost of problem \eqref{problem:aux} is
\begin{equation*}
  \tavg{ f(y) } = \tavg{ f( \lim_{T \rightarrow \infty} \tavg{x^\ast}(T)) } = f( \lim_{T \rightarrow \infty} \tavg{x^\ast}(T)) = f^\optvar.
\end{equation*}
This implies that $\hat{f}^\optvar \leq f^\optvar$.

Alternatively, let $\prtc{ x^\ast(t), y^\ast(t) }_{t=0}^\infty$ be an optimal solution to problem \eqref{problem:aux} such that:
\begin{align*}
  \tavg{ f( y^\ast ) } & = \hat{f}^\optvar \\
  \tavg{ g_j( y^\ast ) } & \leq 0 && j \in \prtc{1, \dotsc, J} \\
  \lim_{T \rightarrow \infty} \tavg{x^\ast}(T) & = \lim_{T \rightarrow \infty} \tavg{y^\ast}(T) \\
  (x^\ast(t), y^\ast(t)) & \in \set{X}_{\random(t)} \times \set{Y} && t \in \prtc{0, 1, 2, \dotsc }. \\
\end{align*}
Consider a solution $\prtc{ x(t) }_{t=0}^\infty$ to problem \eqref{problem:sto} as follows:
\begin{equation*}
  x(t) = x^\ast(t) \quad\quad t \in \prtc{0, 1, 2, \dotsc}.
\end{equation*}
It is easy to see that the solution satisfies the last constraint of problem \eqref{problem:sto}. For the first constraint, the convexity of $g_j$ implies that for $j \in \prtc{1, \dotsc, J}$
\begin{equation*}
  g_j( \lim_{T \rightarrow \infty} \tavg{x}(T) ) = g_j(\lim_{T \rightarrow \infty} \tavg{y^\ast}(T)) \leq \tavg{ g_j(y^\ast) } \leq 0.
\end{equation*}
Hence, this solution is feasible.  The objective cost of problem \eqref{problem:sto} follows from the convexity of $f$ that
\begin{equation*}
  f( \lim_{T \rightarrow \infty} \tavg{x}(T) ) = f( \lim_{T \rightarrow \infty} \tavg{y^\ast}(T) ) \leq \tavg{ f(y^\ast) } = \hat{f}^\optvar.
\end{equation*}
This implies that $f^\optvar \leq \hat{f}^\optvar$.  Thus, combining the above results, we have that $\hat{f}^\optvar = f^\optvar$.
\end{IEEEproof}
\end{theorem}

\subsection{Lyapunov optimization}
The auxiliary problem \eqref{problem:aux} can be solved by the Lyapunov optimization technique \cite{Neely:SNO}.  Define $W_j(t)$ and $Z_i(t)$ to be \emph{virtual queues} of constraints $\expect{ \tavg{ g_j(y) } } \leq 0$ and $\tavg{x}_i  =  \tavg{y}_i$ with update dynamics:
\begin{align}
  W_j(t+1) & = \prts{ W_j(t) + g_j(y(t)) }_+ && j \in \prtc{1, \dotsc, J} \label{eq:W_update} \\
  Z_i(t+1) & = Z_i(t) + x_i(t) - y_i(t) && i \in \prtc{1, \dotsc, I}, \label{eq:Z_update}
\end{align}
where operator $\prts{ \cdot }_+$ is the projection to a corresponding non-negative orthant.  

For ease of notations, let $W(t) \defeq \prtr{ W_1(t), \dotsc, W_J(t) }$, $Z(t) \defeq \prtr{ Z_1(t), \dotsc, Z_I(t) }$, and $g(y) \defeq \prtr{ g_1(y), \dotsc, g_J(y) }$ respectively be the vectors of virtual queues $W_j(t)$, $Z_i(t)$, and functions $g_j(y)$.

Define Lyapunov function \eqref{eq:lya_function} and Lyapunov drift \eqref{eq:lya_drift} as
\begin{align}
  L(t) & \defeq \frac{1}{2} \prts{ \norm{ W(t) }^2 + \norm{ Z(t) }^2 } \label{eq:lya_function} \\
  \drift(t) & \defeq L(t+1) - L(t). \label{eq:lya_drift}
\end{align}
Let notation $A^\tr$ denote the transpose of vector $A$.  Define $C \defeq \sup_{x \in \set{C}, y \in \set{Y}} \prts{ \norm{g(y)}^2 + \norm{x - y}^2 }/2$.

\begin{lemma}
  \label{lem:drift_bound}
  For every $t \in \prtc{0, 1, 2, \dotsc }$, the Lyapunov drift is upper bounded by
\begin{equation}
  \drift(t) \leq C + W(t)^\tr g(y(t)) + Z(t)^\tr \prts{ x(t) - y(t) }.
\end{equation}

\begin{IEEEproof}
From dynamics \eqref{eq:W_update} and \eqref{eq:Z_update}, it follows from the non-expansive projection \cite{Bertsekas:Convex} that
\begin{align*}
  & \norm{W(t+1)}^2 \leq \norm{W(t) + g(y(t))}^2 \\
  & \quad = \norm{W(t)}^2 + \norm{g(y(t))}^2 + 2W(t)^\tr g(y(t))\\
  & \norm{Z(t+1)}^2 = \norm{Z(t) + x(t) - y(t)}^2 \\
  & \quad = \norm{Z(t)}^2 + \norm{x(t) - y(t)}^2 + 2Z(t)^\tr \prts{x(t) - y(t)}.
\end{align*}
Adding the above relations and using definitions \eqref{eq:lya_function} and \eqref{eq:lya_drift} yields
\begin{multline*}
  2 \drift(t) \leq \norm{g(y(t))}^2 + 2W(t)^\tr g(y(t)) \\+ \norm{x(t) - y(t)}^2 + 2Z(t)^\tr \prts{x(t) - y(t)}
\end{multline*}
Using the definition of $C$ proves the lemma.
\end{IEEEproof}
\end{lemma}

Let $V > 0$ be any positive real number representing a parameter of an algorithm solving problem \eqref{problem:aux}.  The drift-plus-penalty term is defined as $\drift(t) + V f(y(t))$.  Applying Lemma \ref{lem:drift_bound}, the drift-plus-penalty term is bounded for every time $t$ by
\begin{multline}
  \label{eq:dpp_bound}
  \drift(t) + V f(y(t)) \leq C + W(t)^\tr g(y(t)) + Z(t)^\tr \prts{ x(t) - y(t) } \\+ V f(y(t)).
\end{multline}

\subsection{Drift-plus-penalty algorithm}
Let $W^0$ and $Z^0$ be the initial condition of $W(0)$ and $Z(0)$ respectively.  Every time step, the Lyapunov optimization technique observes the current realization of random state $\random(t)$ before choosing decisions $x(t) \in \set{X}_{\random(t)}$ and $y(t) \in \set{Y}$ that minimize the right-hand-side of \eqref{eq:dpp_bound}.  The drift-plus-penalty algorithm is summarized in Algorithm \ref{alg:dpp}.

%Then, it minimizes the right-hand-side of \eqref{eq:dpp_bound} with respect to $x \in \set{X}_{\random(t)}$ and $y \in \set{Y}$, and update the virtual queues as dynamics \eqref{eq:W_update} and \eqref{eq:Z_update}.  Let $W^0$ and $Z^0$ be the initial condition of $W(0)$ and $Z(0)$ respectively.  The drift-plus-penalty algorithm is summarized in Algorithm \ref{alg:dpp}.

\begin{algorithm}
  \DontPrintSemicolon
  Initialize $V, W(0) = W^0, Z(0) = Z^0$.\;
  \For{ $t = 0, 1, 2, \dotsc $ } {
    Observe $\random(t)$\;
    $x(t) = \arginf_{x \in \set{X}_{\random(t)}} Z(t)^\tr x$\;
    $y(t) = \arginf_{y \in \set{Y}} [V f(y) + W(t)^\tr g(y) - Z(t)^\tr y]$\;
    $W(t+1) = \prts{ W(t) + g(y(t)) }_+$\;
    $Z(t+1) = Z(t) + x(t) - y(t)$\;
  }
  \caption{Drift-plus-penalty algorithm solving \eqref{problem:aux}.}\label{alg:dpp}
\end{algorithm}

%Define $\result(t) \defeq \prtr{ x(t), y(t), W(t), Z(t) }$ to be the generated results of Algorithm \ref{alg:dpp} at time $t$.  
%Next section shows that a sequence of the generated results exhibits a transient phase and a steady state phase.

\section{Behaviors of Drift-Plus-Penalty Algorithm}
\label{sec:behavior}
Starting from $(W(0), Z(0))$, Algorithm \ref{alg:dpp} reaches the steady state when vector $(W(t),Z(t))$ concentrates around a specific set (defined in Section \ref{sec:embbed}).  The transient phase is the period before this concentration.

\subsection{Embedded Formulation}
\label{sec:embbed}
A convex optimization problem, called \emph{embedded formulation}, is considered.  This idea is inspired by \cite{Longbo:Delay_reduction}.
\begin{align}
  \label{problem:embedded}
  \minimize \quad & f(y) \\
  \subjectto \quad & g_j(y) \leq 0 && j \in \prtc{1, \dotsc, J} \notag\\
  & y = \sum_{\random \in \setrandom} \ssp_\random x^\random \notag \\
  & y,  \in \set{Y},~~ x^\random \in \chull{\set{X}}_\random && \random \in \setrandom. \notag
\end{align}

This formulation has a dual problem, whose properties are used in convergence analysis.  Let $w \in \RealSet_+^J$ and $z \in \RealSet^I$ be the vectors of dual variables associated with the first and second constraints of problem \eqref{problem:embedded}.  The Lagrangian is defined as
\begin{multline*}
  \lagran \prtr{ \prtc{x^\random}_{\random \in \setrandom}, y, w, z } = \\\sum_{\random \in \setrandom} \ssp_\random \prts{ f(y) + w^\tr g(y) + z^\tr (x^\random - y) }.
\end{multline*}
The dual function of problem \eqref{problem:embedded} is
\begin{align}
  d(w, z) & = \inf_{y \in \set{Y},~ x^\random \in \chull{\set{X}}_\random: \forall \random \in \setrandom} \lagran \prtr{ \prtc{x^\random}_{\random \in \setrandom}, y, w, z } \notag\\
%  &\hspace{-2em} = \sum_{\random \in \setrandom} \ssp_\random \inf_{y^\random \in \set{Y},~ x^\random \in \chull{\set{X}}_\random } \prts{ f(y^\random) + w^\tr g(y^\random) + z^\tr (x^\random - y^\random) } \notag\\
  & = \sum_{\random \in \setrandom} \ssp_\random d_\random(w, z), \label{eq:dualfunction}
\end{align}
where $d_{\random}(w, z)$ is defined in \eqref{eq:dualfunction_single} and all of the minimums $y$ take the same value.
\begin{equation}
  \label{eq:dualfunction_single}
  d_\random(w, z) \defeq \inf_{y \in \set{Y},~ x \in \chull{\set{X}}_\random } \prts{ f(y) + w^\tr g(y) + z^\tr (x - y) }.
\end{equation}
Define the solution to the infimum in \eqref{eq:dualfunction_single} as
\begin{align}
  y^\ast(w, z) & \defeq \arginf_{y \in \set{Y}} f(y) + w^\tr g(y) - z^\tr y, \label{eq:yrant}\\
  x^{\ast\random}(z) & \defeq \arginf_{x \in \chull{\set{X}}_\random} z^\tr x. \label{eq:xrant}
\end{align}
%The subgradient of \eqref{eq:dualfunction_single} is
%\begin{equation*}
%  \partial d_\random(w, z) = (g(y^\ast(w,z)), x^{\ast\random}(z) - y^\ast(w,z)).
%\end{equation*}
%Further, $d_{\random}(\dvar)$ in \eqref{eq:dualfunction_single} is concave and satisfies
%\begin{equation}
%  \label{eq:d_concave}
%  d_{\random}( \dvar^\ast ) \leq d_{\random}( \dvar ) + \partial d_{\random}( \dvar )^\tr \prts{ \dvar^\ast - \dvar }.
%\end{equation}

Finally, the dual problem of formulation \eqref{problem:embedded} is
\begin{align}
  \label{problem:dual}
  \maximize \quad & d(w, z) \\
  \subjectto \quad & (w, z) \in \RealSet_+^J \times \RealSet^I. \notag
\end{align}
Problem \eqref{problem:dual} has an optimal solution that may not be unique.  A set of these optimal solutions, which are vectors of Lagrange multipliers, can be used to analyze the expected transient time.  However, to simplify the proofs and notations, the uniqueness assumption is assumed.  Let $\dvar \defeq (w, z)$ denote a concatenation vector of $w$ and $z$.

\begin{assumption}
  \label{ass:uniqueness}
Dual problem \eqref{problem:dual} has a unique vector of Lagrange multipliers denoted by $\dvar^\ast \defeq (w^\ast, z^\ast)$.
\end{assumption}

This assumption is assumed throughout Section \ref{sec:poly} and Section \ref{sec:nonpoly}.  Note that this is a mild assumption when practical systems are considered, e.g., \cite{Longbo:Delay_reduction,Eryilmaz:QLB}.  Furthermore, simulation results in Section \ref{sec:simulation} evince that this assumption may not be needed.

To prove the main result of this section, a useful property of $d_\random(w, z)$ is derived.  Define $h(x,y) \defeq (g(y), x-y)$.
\begin{lemma}
  For any $\dvar = (w, z) \in \RealSet_+^J \times \RealSet^I$ and $\random \in \setrandom$, it holds that
\begin{equation}
  \label{eq:d_concave}
  d_{\random}( \dvar^\ast ) \leq d_{\random}( \dvar ) + h(x^{\ast\random}(z), y^\ast(w, z))^\tr \prts{ \dvar^\ast - \dvar }.
\end{equation}
\begin{IEEEproof}
From \eqref{eq:dualfunction_single}, it follows, for any $\dvar = (w,z) \in \RealSet_+^J \times \RealSet^I$ and $(x,y) \in \chull{\set{X}}_\random \times \set{Y}$, that
\begin{align*}
  d_\random(\dvar^\ast) & \leq f(y) + h(x, y)^\tr \dvar^\ast  \\
  & = f(y) + h(x, y)^\tr \dvar  + h(x, y)^\tr \prts{ \dvar^\ast - \dvar }
\end{align*}
Setting $(x, y) = \prtr{x^{\ast\random}(w, z), y^{\ast\random}(z)}$, as defined in \eqref{eq:yrant} and \eqref{eq:xrant}, and using \eqref{eq:dualfunction_single} proves the lemma.
\end{IEEEproof}
\end{lemma}

The following lemma ties the virtual queues of Algorithm \ref{alg:dpp} to the Lagrange multipliers.  Given the generated result of Algorithm \ref{alg:dpp}, define $Q(t) \defeq (W(t), Z(t))$ as a concatenation of vectors $W(t)$ and $Z(t)$.

\begin{lemma}
  \label{lem:onestep}
  The following holds for every $t \in \prtc{0, 1, 2, \dotsc}$:
\begin{multline*}
  \expect{ \norm{ Q(t+1) - V \dvar^\ast }^2 | Q(t) } \\\leq \norm{ Q(t) - V \dvar^\ast }^2 + 2C + 2V \prts{ d( Q(t)/V ) - d( \dvar^\ast ) }.
\end{multline*}

\begin{IEEEproof}
The non-expansive projection \cite{Bertsekas:Convex} implies that
\begin{align}
  & \norm{ Q(t+1) - V \dvar^\ast }^2 \leq \norm{ Q(t) + h(x(t), y(t)) - V\dvar^\ast }^2 \notag\\
  &\hspace{-0.7em} \quad = \norm{ Q(t) - V\dvar^\ast }^2 + \norm{ h(x(t), y(t)) }^2 \notag\\
  &\hspace{-0.7em} \quad\quad\quad + 2 h(x(t), y(t))^\tr \prts{ Q(t) - V\dvar^\ast } \notag\\
  &\hspace{-0.7em} \quad \leq \norm{ Q(t) - V\dvar^\ast }^2 + 2C + 2 h(x(t), y(t))^\tr \prts{ Q(t) - V\dvar^\ast } \label{lem_os:1}
\end{align}

From \eqref{eq:dualfunction_single}, when $\dvar = Q(t)/V$, we have 
\begin{multline*}
d_{\random(t)}( Q(t)/V ) = \\\inf_{y \in \set{Y},~ x \in \chull{\set{X}}_{\random(t)} } \prts{ f(y) + \frac{W(t)}{V}^\tr g(y) + \frac{Z(t)}{V}^\tr (x - y) },
\end{multline*}
so $y^\ast(W(t)/V, Z(t)/V) = y(t)$ and $x^{\ast\random(t)}(Z(t)/V) = x(t)$ where $(y^\ast(W(t)/V,Z(t)/V), x^{\ast\random(t)}(Z(t)/V))$ is defined in \eqref{eq:yrant} and \eqref{eq:xrant}, and $(y(t), x(t))$ is the decision from Algorithm \ref{alg:dpp}.  Therefore, property \eqref{eq:d_concave} implies that
\begin{equation*}
  h(x(t), y(t))^\tr \prts{ Q(t) - V\dvar^\ast } \leq V \prts{ d_{\random(t)}( Q(t)/V ) - d_{\random(t)}( \dvar^\ast ) }.
\end{equation*}
Applying the above inequality on the last term of \eqref{lem_os:1} gives
\begin{multline*}
  \norm{ Q(t+1) - V \dvar^\ast }^2 \leq \norm{ Q(t) - V\dvar^\ast }^2 + 2C \\+ 2V \prts{ d_{\random(t)}( Q(t)/V ) - d_{\random(t)}( \dvar^\ast ) }.
\end{multline*}
Taking a conditional expectation given $Q(t)$ proves the lemma.
\begin{multline*}
  \expect{\norm{ Q(t+1) - V \dvar^\ast }^2 | Q(t) } \leq \norm{ Q(t) - V\dvar^\ast }^2 + 2C \\+ 2V \sum_{\random \in \setrandom} \ssp_\random \prts{ d_{\random}( Q(t)/V ) - d_{\random}( \dvar^\ast ) }.
\end{multline*}
\end{IEEEproof}
\end{lemma}

The analysis of transient and steady state phases in Sections \ref{sec:poly} and \ref{sec:nonpoly} will utilize Lemma \ref{lem:onestep}.  The convergence results in the steady state require the following results.

\subsection{$T$-slot convergence}
For any positive integer $T$ and any starting time $t_0$, define the $T$-slot average starting at $t_0$ as
\begin{equation*}
  \tavg{x}(t_0, T) \defeq \frac{1}{T} \sum_{t = t_0}^{t_0 + T -1} x(t).
\end{equation*}
This average leads to the following convergence bounds.

\begin{theorem}
  \label{thm:gen_convergence}
  Let $\prtc{ Q(t) }_{t=0}^\infty$ be a sequence generated by Algorithm \ref{alg:dpp}.  For any positive integer $T$ and any starting time $t_0$, the objective cost converges as
\begin{multline}
  \label{eq:obj_gen}
  \expect{ f( \tavg{x}(t_0, T) ) } - f^\optvar \leq  \frac{M_f}{T} \expect{ \norm{ Z(t_0 + T) - Z(t_0) } } \\ + \frac{1}{2TV} \expect{ \norm{Q(t_0)}^2 - \norm{Q(t_0 + T)}^2 } + \frac{C}{V},
\end{multline}
and the constraint violation for every $j \in \prtc{1, \dotsc, J}$ is
\begin{multline}
  \label{eq:con_gen}
  \expect{ g_j(\tavg{x}(t_0, T)) } \leq \frac{1}{T} \expect{ W_j(t_0 + T) - W_j(t_0) } \\+ \frac{M_{gj}}{T} \expect{ \norm{ Z(t_0 + T) - Z(t_0) } }.
\end{multline}
\begin{IEEEproof}
The proof is in Appendix.
\end{IEEEproof}
\end{theorem}

To interprete Theorem \ref{thm:gen_convergence}, the following concentration bound is provided.  It is proven in \cite{Neely:1Queue}.

\subsection{Concentration bound}
%Let $\indicator{A}$ be an indicator function whose value is $1$ if $A$ is true and $0$ otherwise.
\begin{theorem}
  \label{thm:concentration}
Let $K(t)$ be a real random process over $t \in \prtc{0, 1, 2, \dotsc}$ satisfying
\begin{align*}
  \abs{ K(t+1) - K(t) } & \leq \delta \quad\quad \text{and}\\
  \expect{ K(t+1) - K(t) | K(t) } & \leq \left\{
    \begin{array}{ll}
      \delta & , K(t) < K \\
      -\beta & , K(t) \geq K,
    \end{array}
    \right.
\end{align*}
for some positive real-valued $\delta, K$, and $0 < \beta \leq \delta$.

Suppose $K(0) = k_0$ (with probability $1$) for some $k_0 \in \RealSet$.  Then for every time $t \in \prtc{0, 1, 2, \dotsc}$, the following holds:
\begin{equation*}
  \expect{ e^{rK(t)} } \leq D + \prtr{ e^{rk_0} - D } \rho^t
\end{equation*}
where $0 < \rho < 1$ and constants $r, \rho$, and $D$ are:
\begin{align*}
  r & \defeq \frac{\beta}{\prtr{\delta^2 + \delta\beta/3}}, & \rho & \defeq 1 - \frac{r\beta}{2} \\
  D & \defeq \frac{ \prtr{ e^{r\delta} - \rho } e^{rK} }{ 1 - \rho }.
\end{align*}
\end{theorem}

In this paper, random process $K(t)$ is defined to be the distance between $Q(t)$ and the vector of Lagrange multipliers as $K(t) \defeq \norm{ Q(t) - V\dvar^\ast }$ for every $t \in \prtc{0, 1, 2, \dotsc}$.

\begin{lemma}
  \label{lem:k}
  It holds for every $t \in \prtc{0, 1, 2, \dotsc}$ that
\begin{align*}
  \abs{ K(t+1) - K(t) } & \leq \sqrt{2C} \\
  \expect{ K(t+1) - K(t) | K(t) } & \leq \sqrt{2C}.
\end{align*}

\begin{IEEEproof}
The first part is proven in two cases.  

i) If $K(t+1) \geq K(t)$, the non-expansive projection implies
\begin{align*}
  & \abs{K(t+1) - K(t)} = K(t+1) - K(t) \\
  & \quad \leq \norm{ Q(t) + h(x(t), y(t)) - V\dvar^\ast } - \norm{ Q(t) - V\dvar^\ast } \\
  & \quad \leq \norm{ h(x(t), y(t)) } \leq \sqrt{2C}.
\end{align*}

ii) If $K(t+1) < K(t)$, then 
\begin{align*}
  & \abs{K(t+1) - K(t)} = K(t) - K(t+1) \\
  & \quad \leq \norm{ Q(t) - Q(t+1) } + \norm{ Q(t+1) - V\dvar^\ast } - K(t+1) \\
  & \quad \leq \norm{ h(x(t), y(t)) } \leq \sqrt{2C}.
\end{align*}
Therefore, $\abs{ K(t+1) - K(t) } \leq \sqrt{2C}$.  Using $K(t+1) - K(t) \leq \abs{ K(t+1) - K(t) }$ proves the second part.
\end{IEEEproof}
\end{lemma}

Lemma \ref{lem:k} prepares $K(t)$ for Theorem \ref{thm:concentration}.  The only constants left to be specified are $\beta$ and $K$, which depend on properties of dual function \eqref{eq:dualfunction}.

\section{Locally-polyhedral dual function}
\label{sec:poly}
\begin{figure}
  \centering
  \includegraphics[scale=0.78]{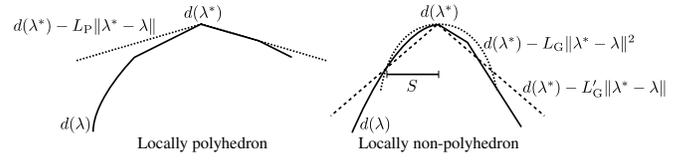}
  \caption{Illustrations of locally-polyhedral and locally-non-polyhedral dual functions}
  \label{fig:local_property}
\end{figure}

This section analyzes the expected transient time and a convergence result in the steady state.  Dual function \eqref{eq:dualfunction} in this section is assumed to satisfy a locally-polyhedral property, introduced in \cite{Longbo:Delay_reduction}.  This property is illustrated in Figure \ref{fig:local_property}.  It holds when $f$ and each $g_j$ for $j \in \prtc{1, \dotsc, J}$ are either linear or piece-wise linear.

\begin{assumption}
  \label{ass:polyhedral}
  Let $\dvar^\ast$ be the unique Lagrange multiplier vector.  There exists $L_\poly > 0$ such that dual function \eqref{eq:dualfunction} satisfies, for any $\dvar \in \RealSet_+^J \times \RealSet^I$,
\begin{equation}
  \label{eq:polyhedral}
  d(\dvar^\ast) \geq d(\dvar) + L_\poly \norm{ \dvar - \dvar^\ast }.
\end{equation}
\end{assumption}

Note that, by concavity of the dual function, if inequality \eqref{eq:polyhedral} holds locally about $\dvar^\ast$, it must also hold globally.

\subsection{Expected transient time}
The progress of $Q(t)$ at each step can be analyzed.  Define 
\begin{equation*}
  B_\poly \defeq \max\prts{ \frac{L_\poly}{2}, \frac{2C}{L_\poly} }.
\end{equation*}

\begin{lemma}
  \label{lem:poly_onestep}
Under Assumptions \ref{ass:uniqueness} and \ref{ass:polyhedral}, whenever $\norm{ Q(t) - V\dvar^\ast } \geq B_\poly$, the following holds
\begin{equation*}
  \label{eq:poly_onestep}
  \expect{ \norm{ Q(t+1) - V\dvar^\ast } | Q(t) } \leq \norm{ Q(t) - V\dvar^\ast } - \frac{L_\poly}{2}.
\end{equation*}
\begin{IEEEproof}
If condition
\begin{equation}
  \label{lem_po:1}
  2C + 2V\prts{ d( Q(t)/V ) - d( \dvar^\ast ) } \leq -2\eta \norm{ Q(t) - V\dvar^\ast } + \eta^2
\end{equation}
is true, Lemma \ref{lem:onestep} implies that
\begin{align*}
  & \expect{ \norm{ Q(t+1) - V \dvar^\ast }^2 | Q(t) } \\
  & \quad \leq \norm{ Q(t) - V \dvar^\ast }^2 - 2\eta \norm{ Q(t) - V\dvar^\ast } + \eta^2 \\
  & \quad = \prts{ \norm{ Q(t) - V \dvar^\ast } - \eta }^2.
\end{align*}
Applying Jensen's inequality \cite{Bertsekas:Convex} on the left-hand-side yields
\begin{equation*}
  \expect{ \norm{ Q(t+1) - V \dvar^\ast } | Q(t) }^2 \leq \prts{ \norm{ Q(t) - V \dvar^\ast } - \eta }^2.
\end{equation*}
When $\norm{ Q(t) - V\dvar^\ast } \geq \eta$, it follows that
\begin{equation}
  \label{lem_po:3}
  \expect{ \norm{ Q(t+1) - V \dvar^\ast } | Q(t) } \leq \norm{ Q(t) - V \dvar^\ast } - \eta.
\end{equation}

However, condition \eqref{lem_po:1} holds when
\begin{equation}
  \label{lem_po:2}
  -2VL_\poly \norm{ Q(t)/V - \dvar^\ast } \leq -2\eta \norm{ Q(t) - V\dvar^\ast } - 2C,
\end{equation}
because Assumption \ref{ass:polyhedral} implies that
\begin{equation*}
  2V \prts{ d(Q(t)/V) - d(\dvar^\ast) } \leq -2VL_\poly \norm{ Q(t)/V - \dvar^\ast }.
\end{equation*}
Therefore, condition \eqref{lem_po:1} holds when condition \eqref{lem_po:2} holds.  This requires that $\norm{ Q(t) - V\dvar^\ast } \geq C/(L_\poly - \eta)$ when $\eta < L_\poly$.  Thus, inequality \eqref{lem_po:3} holds when $\norm{ Q(t) - V\dvar^\ast } \geq \max\prts{ \eta, C/(L_\poly - \eta) }$.  Choosing $\eta = L_\poly/2$ proves the lemma.
\end{IEEEproof}
\end{lemma}

Lemma \ref{lem:poly_onestep} implies that $Q(t)$ proceeds closer to $V\dvar^\ast$ in the next step when the distance between them is at least $B_\poly$.  This implication means that $Q(t)$ concentrates around $V\dvar^\ast$ in the steady state.  The transient time $T_\poly$ of Algorithm \ref{alg:dpp} under the locally-polyhedral assumption is defined as
\begin{equation*}
  T_\poly \defeq \arginf_{t \geq 0} \prtc{ \norm{ Q(t) - V\dvar^\ast } < B_\poly }.
\end{equation*}
Also, let $\indicator{A}$ be an indicator function whose value is $1$ if $A$ is true and $0$ otherwise.

\begin{lemma}
  \label{lem:transient_poly}
  Given the initial condition $Q^0 \defeq (W^0, Z^0)$, the expected transient time is at most
\begin{equation*}
  \expect{ T_\poly } \leq \frac{2 \norm{ Q^0 - V\dvar^\ast }}{L_\poly}.
\end{equation*}
\begin{IEEEproof}
Recall that $K(t) \defeq \norm{ Q(t) - V\dvar^\ast }$.  Define random variable $S(t)$ as $S(t) = K(t+1) - K(t)$, and $E_n = \min\prts{ n, T_\poly-1 }$.  It follows that
\begin{align*}
  \sum_{t=0}^{n} S(t) \indicator{t \leq T_\poly - 1} &= \sum_{t=0}^{E_n} S(t) \\
  &= K(E_n+1) - K(0) \geq -K(0)
\end{align*}
where $K(0) = \norm{Q^0 - V\dvar^\ast}$ is a given constant (with probability 1).  Taking an expectation of both sides gives:
\begin{align*}
  \sum_{t=0}^n \expect{ S(t) | t \leq T_\poly - 1 } \prob{t \leq T_\poly - 1} \geq -K(0)
\end{align*}

However, when $t \leq T_\poly - 1$, we know that $\norm{Q(t) - V\dvar^\ast} \geq B_\poly$ from the definition of $T_\poly$.  Then, Lemma \ref{lem:poly_onestep} implies that $\expect{ S(t) | t \leq T_\poly - 1 } \leq -L_\poly/2$, and
\begin{equation*}
  -\frac{L_\poly}{2} \sum_{t=0}^n \prob{t \leq T_\poly - 1} \geq -K(0).
\end{equation*}
Taking a limit as $n \rightarrow \infty$ gives:
\begin{equation*}
  -\frac{L_\poly}{2} \expect{T_\poly} \geq -K(0),
\end{equation*}
where $\expect{T_\poly} = \sum_{t=0}^\infty \prob{ T_\poly > t }$.  This proves the lemma.  Note that this argument is similar to a drift lemma given in Proposition 5.3 in \cite{Soren:prob}.
\end{IEEEproof}
\end{lemma}
 
Lemma \ref{lem:transient_poly} implies that the expected transient time under the locally-polyhedral assumption is $O(V)$.

\subsection{Convergence time in a steady state}
Once Algorithm \ref{alg:dpp} enters the steady state, the following property holds.

\begin{lemma}
  \label{lem:distance_poly}
  Under Assumptions \ref{ass:uniqueness} and \ref{ass:polyhedral}, for any time $t \geq T_\poly$, the following holds 
\begin{align}
  \expect{ \norm{ Q(t) - V\dvar^\ast } } & \leq U_\poly \label{eq:ex_poly}\\
  \expect{ \norm{ Q(t) - V\dvar^\ast }^2 } & \leq U_\poly' \label{eq:ex2_poly}
\end{align}
where constants $U_\poly, U_\poly', r_\poly, \rho_\poly,$ and $D_\poly$ are:
\begin{align*}
  U_\poly &\defeq \frac{ \log\prtr{D_\poly + e^{r_\poly B_\poly}} }{ r_\poly }, & U_\poly' &\defeq \frac{ 2\prtr{ D_\poly + e^{r_\poly B_\poly} } }{ r_\poly^2 }, \\
  r_\poly &\defeq \frac{3L_\poly}{12C + L_\poly\sqrt{2C}}, & \rho_\poly &\defeq 1 - \frac{r_\poly L_\poly}{4}, \\
  D_\poly &\defeq \frac{e^{r_\poly B_\poly}\prtr{ e^{r_\poly \sqrt{2C}} - \rho_\poly}}{ 1 - \rho_\poly }.
\end{align*}

\begin{IEEEproof}
Recall that $K(t) \defeq \norm{ Q(t) - V\dvar^\ast }$.  From Lemmas \ref{lem:k} and \ref{lem:poly_onestep}, constants in Theorem \ref{thm:concentration} are $\delta = \sqrt{2C}, K = B_\poly$, and $\beta = L_\poly/2$.  By the definition of $T_\poly$, $k_0 = K(T_\poly) < K$, and Theorem \ref{thm:concentration} implies, for any $t \geq T_\poly$, that
\begin{equation}
  \label{lem_distance:1}
  \expect{ e^{rK(t)} } \leq D + \prtr{ e^{rk_0} - D }\rho^\prtr{t - T_\poly} \leq D + e^{rk_0}.
\end{equation}
Jensen's inequality implies that $e^{r\expect{K(t)}} \leq \expect{ e^{rK(t)} }$, and we have $e^{r\expect{K(t)}} \leq D + e^{rk_0}$.  Taking logarithm and dividing by $r$ proves \eqref{eq:ex_poly}.

Chernoff bound \cite{Ross:Stochastic} implies that, for any $m \in \RealSet_+$,
\begin{equation}
  \label{lem_distance:2}
  \prob{ K(t) \geq m } \leq e^{-rm} \expect{ e^{r K(t) } } \leq \prtr{D + e^{rk_0}} e^{-rm},
\end{equation}
where the last inequality uses \eqref{lem_distance:1}.

Since $K(t)^2$ is always non-negative, it can be shown that $\expect{ K(t)^2 } = 2 \int_0^\infty m \prob{ K(t) \geq m } dm$ by the integration by parts.  Using \eqref{lem_distance:2}, we have
\begin{align*}
  \expect{ K(t)^2 } &\leq 2\prtr{D + e^{rk_0}} \int_0^\infty m e^{-rm} dm.
\end{align*}
Performing the integration by parts proves \eqref{eq:ex2_poly}.
\end{IEEEproof}
\end{lemma}

The above lemma implies that, in the steady state, the expected distance and square distance between $Q(t)$ and the vector of Lagrange multipliers are constant.  This phenomenon leads to an improved convergence time when the average is performed in the steady state.  A useful result is derived before the main theorem.

\begin{lemma}
  \label{lem:Qsquare}
  For any times $t_1$ and $t_2$, it holds that
\begin{multline*}
  \expect{ \norm{Q(t_1)}^2 - \norm{Q(t_2)}^2 } \leq \expect{ \norm{ Q(t_1) - V\dvar^\ast }^2 } \\+ 2\norm{ V\dvar^\ast } \expect{ \norm{ Q(t_1) - V\dvar^\ast } + \norm{ Q(t_2) - V\dvar^\ast } }.
\end{multline*}
\begin{IEEEproof}
It holds for any $Q \in \RealSet_+^J \times \RealSet^I$ that
\begin{equation*}
  \norm{Q}^2 = \norm{Q-V\dvar^\ast}^2 + \norm{V\dvar^\ast}^2 + 2(Q-V\dvar^\ast)^\tr(V\dvar^\ast).
\end{equation*}
Using the above equality with $Q_1, Q_2 \in \RealSet_+^J \times \RealSet^I$ leads to
\begin{align*}
  & \norm{Q_1}^2 - \norm{Q_2}^2 \\
  & \quad \leq \norm{Q_1-V\dvar^\ast}^2 - \norm{Q_2 - V\dvar^\ast}^2 + 2(Q_1 - Q_2)^\tr(V\dvar^\ast) \\
  & \quad \leq \norm{Q_1-V\dvar^\ast}^2 + 2\norm{Q_1 - Q_2} \norm{V\dvar^\ast} \\
  & \quad \leq \norm{Q_1-V\dvar^\ast}^2 + 2\norm{V\dvar^\ast}\prts{ \norm{Q_1-V\dvar^\ast} + \norm{Q_2-V\dvar^\ast} }.
\end{align*}
Taking an expectation proves the lemma.
\end{IEEEproof}
\end{lemma}

Finally, the convergence in the steady state is analyzed.
\begin{theorem}
  \label{thm:poly}
  For any time $t_0 \geq T_\poly$ and positive integer $T$, the objective cost converges as
\begin{equation}
  \label{eq:obj_poly}
  \expect{ f( \tavg{x}(t_0, T) ) } - f^\optvar \leq \frac{2M_f U_\poly}{T} + \frac{U_\poly' + 4VU_\poly\norm{\dvar^\ast}}{2TV} + \frac{C}{V}
\end{equation}
and the constraint violation is upper bounded by
\begin{equation}
  \label{eq:con_poly}
  \expect{ g_j(\tavg{x}(t_0, T)) } \leq \frac{2U_\poly}{T} + \frac{2M_{gj} U_\poly}{T}.
\end{equation}
\begin{IEEEproof}
From Theorem \ref{thm:gen_convergence}, the objective cost converges as \eqref{eq:obj_gen}.  Since $T_\poly \leq t_0 < t_0+T$, we use results in Lemma \ref{lem:distance_poly} to upper bound $\expect{ \norm{ Q(t) - V\dvar^\ast } }$ and $\expect{ \norm{ Q(t) - V\dvar^\ast }^2 }$ for $t_0$ and $t_0 + T$.
%\begin{multline}
%  \label{thm-poly:1}
%  \expect{ f( \tavg{x}(t_0, T) ) } - f^\optvar \leq \frac{M_f}{T} \expect{ \norm{ Z(t_0 + T) - Z(t_0) } } \\ + \frac{1}{2TV} \expect{ \norm{Q(t_0)}^2 - \norm{Q(t_0 + T)}^2 } + \frac{C}{V}.
%\end{multline}
Terms in the right-hand-side of \eqref{eq:obj_gen} are bounded by
\begin{align}
  & \expect{ \norm{ Z(t_0 + T) - Z(t_0) } } \leq \expect{ \norm{ Q(t_0 + T) - Q(t_0) } }  \notag\\
  & \quad \leq \expect{ K(t_0 + T) + K(t_0) } \leq 2U_\poly. \label{thm-poly:2}
\end{align}
Lemma \ref{lem:Qsquare} implies that
\begin{align}
  &\expect{ \norm{Q(t_0)}^2 - \norm{Q(t_0 + T)}^2 } \notag\\
  & \quad \leq \expect{ K(t_0)^2 + 2\norm{ V\dvar^\ast }\prts{ K(t_0) + K(t_0+T) } } \notag\\
  & \quad \leq U_\poly' + 4VU_\poly \norm{\dvar^\ast}. \label{thm-poly:3}
\end{align}
Substituting bounds \eqref{thm-poly:2} and \eqref{thm-poly:3} into \eqref{eq:obj_gen} proves \eqref{eq:obj_poly}.

The constraint violation converges as \eqref{eq:con_gen} where $T_\poly \leq t_0 < t_0 + T$, and Lemma \ref{lem:distance_poly} can be utilized.
%\begin{multline}
%  \label{thm-poly:4}
%  \expect{ g_j(\tavg{x}(t_0, T)) } \leq \frac{1}{T} \expect{ W_j(t_0 + T) - W_j(t_0) } \\+ \frac{M_{gj}}{T} \expect{ \norm{ Z(t_0+T) - Z(t_0) } }.
%\end{multline}
The last term in the right-hand-side of \eqref{eq:con_gen} is bounded in \eqref{thm-poly:2}.  The first term is bounded by
\begin{align*}
  & \expect{ W_j(t_0 + T ) - W_j(t_0) } \\
  & \quad \leq \expect{ \abs{ W_j(t_0 + T ) - V\dvar^\ast } + \abs{ W_j(t_0) - V\dvar^\ast } } \\ 
  & \quad \leq \expect{ K(t_0 + T) + K(t_0) } \leq 2U_\poly.
\end{align*}
Substituting the above bound and \eqref{thm-poly:2} into \eqref{eq:con_gen} proves \eqref{eq:con_poly}.
\end{IEEEproof}
\end{theorem}

The implication of Theorem \ref{thm:poly} is as follows.  When the average starts in the steady state, the deviation from the optimal cost is $O(1/T + 1/V)$, and the constraint violation is bounded by $O(1/T)$.  By setting $V= 1/\epsilon$, the convergence time is $O(1/\epsilon)$.  
%Further, since the expected transient time is $O(1/\epsilon)$, if the average is performed shortly after the transient time, the total time requires to achieve $\epsilon$-optimality is $O(1/\epsilon)$.

\section{Locally-non-polyhedral dual function}
\label{sec:nonpoly}
The dual function \eqref{eq:dualfunction} in Section \ref{sec:nonpoly} is assumed to satisfy a locally-non-polyhedral property, modified from \cite{Longbo:Delay_reduction}.  This property is illustrated in Figure \ref{fig:local_property}.

\begin{assumption}
  \label{ass:nonpolyhedral}
  Let $\dvar^\ast$ be the unique Lagrange multiplier vector.  The following holds:

i) There exist $S > 0$ and $L_\gen > 0$ such that, whenever $\dvar \in \RealSet_+^J \times \RealSet^I$ and $\norm{ \dvar - \dvar^\ast } \leq S$, dual function \eqref{eq:dualfunction} satisfies 
\begin{equation*}
  \label{eq:nonpolyhedral}
  d(\dvar^\ast) \geq d(\dvar) + L_\gen \norm{ \dvar - \dvar^\ast }^2.
\end{equation*}

ii) When $\dvar \in \RealSet_+^J \times \RealSet^I$ and $\norm{ \dvar - \dvar^\ast } > S$, there exist $L_\gen' > 0$ such that dual function \eqref{eq:dualfunction} satisfies
\begin{equation*}
  d(\dvar^\ast) \geq d(\dvar) + L_\gen' \norm{ \dvar - \dvar^\ast }.
\end{equation*}
\end{assumption}

The progress of $Q(t)$ at each step can be analyzed.  Define
\begin{align*}
  B_\gen(V) & \defeq \max\prts{ \frac{1}{\sqrt{V}}, \sqrt{V} \prtr{ \frac{1 + \sqrt{1 + 4L_\gen C}}{2L_\gen} } } \\
  B_\gen' & \defeq \max\prts{ \frac{L_\gen'}{2}, \frac{2C}{L_\gen'} }.
\end{align*}

\begin{lemma}
  \label{lem:gen_onestep}
  When $V$ is sufficiently larger that $B_\gen(V) < SV$ and $B_\gen' \leq SV$, under Assumptions \ref{ass:uniqueness} and \ref{ass:nonpolyhedral}, the following holds
\begin{align}
  & \expect{ \norm{ Q(t+1) - V\dvar^\ast } | Q(t) } - \norm{ Q(t) - V\dvar^\ast } \notag\\
  & \quad \leq \left\{
\begin{array}{ll}
  - \frac{1}{\sqrt{V}} & \quad ,\text{if}~ B_\gen(V) \leq \norm{ Q(t) - V\dvar^\ast } \leq SV \\
  - \frac{L_\gen'}{2} & \quad ,\text{if}~ \norm{ Q(t) - V\dvar^\ast } > SV
\end{array}\right.   \label{eq:gen_onestep}
\end{align}
\begin{IEEEproof}
If condition
\begin{equation}
  \label{lem_gen:1}
  2C + 2V\prts{ d( Q(t)/V ) - d( \dvar^\ast ) } \leq -\frac{2}{\sqrt{V}} \norm{ Q(t) - V\dvar^\ast } + \frac{1}{V}
\end{equation}
is true, Lemma \ref{lem:onestep} implies that
\begin{align*}
  & \expect{ \norm{ Q(t+1) - V \dvar^\ast }^2 | Q(t) } \\
  & \quad \leq \norm{ Q(t) - V \dvar^\ast }^2 - \frac{2}{\sqrt{V}} \norm{ Q(t) - V\dvar^\ast } + \frac{1}{V} \\
  & \quad = \prts{ \norm{ Q(t) - V \dvar^\ast } - \frac{1}{\sqrt{V}} }^2.
\end{align*}
Applying Jensen's inequality \cite{Bertsekas:Convex} on the left-hand-side yields
\begin{equation*}
  \expect{ \norm{ Q(t+1) - V \dvar^\ast } | Q(t) }^2 \leq \prts{ \norm{ Q(t) - V \dvar^\ast } - \frac{1}{\sqrt{V}} }^2.
\end{equation*}
When $\norm{ Q(t) - V\dvar^\ast } \geq 1/\sqrt{V}$, it follows that
\begin{equation}
  \label{lem_gen:3}
  \expect{ \norm{ Q(t+1) - V \dvar^\ast } | Q(t) } \leq \norm{ Q(t) - V \dvar^\ast } - \frac{1}{\sqrt{V}}.
\end{equation}

However, condition \eqref{lem_gen:1} holds when
\begin{equation}
  \label{lem_gen:2}
  -2VL_\gen \norm{ Q(t)/V - \dvar^\ast }^2 \leq -\frac{2}{\sqrt{V}} \norm{ Q(t) - V\dvar^\ast } - 2C,
\end{equation}
because Assumption \ref{ass:nonpolyhedral} implies that, when $\norm{ Q(t) - V\dvar^\ast } \leq SV$,
\begin{equation*}
  2V \prts{ d(Q(t)/V) - d(\dvar^\ast) } \leq -2VL_\gen \norm{ Q(t)/V - \dvar^\ast }^2.
\end{equation*}
Therefore, condition \eqref{lem_gen:1} holds when condition \eqref{lem_gen:2} holds.  Condition \eqref{lem_gen:2} requires that
\begin{equation*}
  \norm{ Q(t) - V\dvar^\ast } \geq \frac{\sqrt{V} + \sqrt{V + 4 L_\gen V C }}{2L_\gen}.
\end{equation*}
Thus, inequality \eqref{lem_gen:3} holds when $\norm{ Q(t) - V\dvar^\ast } \geq \max\prts{ \frac{1}{\sqrt{V}}, \frac{\sqrt{V} + \sqrt{V + 4 L_\gen V C }}{2L_\gen} }$.  This proves the first part of \eqref{eq:gen_onestep}.

For the last part of \eqref{eq:gen_onestep}, if condition
\begin{equation}
  \label{lem_npo:1}
  2C + 2V\prts{ d( Q(t)/V ) - d( \dvar^\ast ) } \leq -2\eta \norm{ Q(t) - V\dvar^\ast } + \eta^2
\end{equation}
is true, Lemma \ref{lem:onestep} implies that
\begin{align*}
  & \expect{ \norm{ Q(t+1) - V \dvar^\ast }^2 | Q(t) } \\
  & \quad \leq \norm{ Q(t) - V \dvar^\ast }^2 - 2\eta \norm{ Q(t) - V\dvar^\ast } + \eta^2 \\
  & \quad = \prts{ \norm{ Q(t) - V \dvar^\ast } - \eta }^2.
\end{align*}
Applying Jensen's inequality \cite{Bertsekas:Convex} on the left-hand-side yields
\begin{equation*}
  \expect{ \norm{ Q(t+1) - V \dvar^\ast } | Q(t) }^2 \leq \prts{ \norm{ Q(t) - V \dvar^\ast } - \eta }^2.
\end{equation*}
When $\norm{ Q(t) - V\dvar^\ast } \geq \eta$, it follows that
\begin{equation}
  \label{lem_npo:3}
  \expect{ \norm{ Q(t+1) - V \dvar^\ast } | Q(t) } \leq \norm{ Q(t) - V \dvar^\ast } - \eta.
\end{equation}

However, condition \eqref{lem_npo:1} holds when
\begin{equation}
  \label{lem_npo:2}
  -2VL_\gen \norm{ Q(t)/V - \dvar^\ast } \leq -2\eta \norm{ Q(t) - V\dvar^\ast } - 2C,
\end{equation}
because Assumption \ref{ass:nonpolyhedral} implies that, when $\norm{ Q(t) - V\dvar^\ast } > SV$,
\begin{equation*}
  2V \prts{ d(Q(t)/V) - d(\dvar^\ast) } \leq -2VL_\gen' \norm{ Q(t)/V - \dvar^\ast }.
\end{equation*}
Therefore, condition \eqref{lem_npo:1} holds when condition \eqref{lem_npo:2} holds.  This requires that $\norm{ Q(t) - V\dvar^\ast } \geq C/(L_\gen' - \eta)$ when $\eta < L_\gen'$.  Thus, inequality \eqref{lem_npo:3} holds when $\norm{ Q(t) - V\dvar^\ast } \geq \max\prts{ \eta, C/(L_\gen' - \eta) }$.  Choosing $\eta = L_\gen'/2$ and using the fact that $B_\gen' \leq SV$ proves the last part of \eqref{eq:gen_onestep}.
\end{IEEEproof}
\end{lemma}

Lemma \ref{lem:gen_onestep} implies that $Q(t)$ proceeds closer to $V\dvar^\ast$ in the next step when the distance between them is at least $B_\gen(V)$.  This implication means that $Q(t)$ concentrates around $V\dvar^\ast$ in the steady state.

\subsection{Expected transient time}
From an initial condition, the transient time $T_\gen$ of Algorithm \ref{alg:dpp} under the locally-non-polyhedral assumption is defined as
\begin{equation*}
  T_\gen \defeq \arginf_{t \geq 0} \prtc{ \norm{ Q(t) - V\dvar^\ast } < B_\gen(V) }.
\end{equation*}

\begin{lemma}
  \label{lem:transient_gen2}
  When $V$ is sufficiently larger that $B_\gen(V) < SV$ and $B_\gen' \leq SV$, given the initial condition $Q^0$, the expected transient time is at most
\begin{equation*}
  \expect{ T_\gen } \leq \max\prts{ \sqrt{V}, \frac{2}{L_\gen'} } \norm{ Q^0 - V\dvar^\ast }.
\end{equation*}
\begin{IEEEproof}
Recall that $K(t) \defeq \norm{ Q(t) - V\dvar^\ast }$.  Define random variable $S(t)$ as $S(t) = K(t+1) - K(t)$, and $E_n = \min\prts{ n, T_\gen-1 }$.  It follows that
\begin{align*}
  \sum_{t=0}^{n} S(t) \indicator{t \leq T_\gen - 1} &= \sum_{t=0}^{E_n} S(t) \\
  &= K(E_n +1) - K(0) \geq -K(0)
\end{align*}
where $K(0) = \norm{Q^0 - V\dvar^\ast}$ is a given constant (with probability 1).  Taking an expectation of both sides gives:
\begin{align*}
  \sum_{t=0}^n \expect{ S(t) | t \leq T_\gen - 1 } \prob{t \leq T_\gen - 1} \geq -K(0)
\end{align*}

However, when $t \leq T_\gen - 1$, we know that $\norm{Q(t) - V\dvar^\ast} \geq B_\gen(V)$ from the definition of $T_\gen$.  Then, Lemma \ref{lem:poly_onestep} implies that $\expect{ S(t) | t \leq T_\gen - 1 } \leq -\min\prts{ \frac{1}{\sqrt{V}}, \frac{L_\gen'}{2} }$, and
\begin{equation*}
  -\min\prts{ \frac{1}{\sqrt{V}}, \frac{L_\gen'}{2} } \sum_{t=0}^n \prob{t \leq T_\gen - 1} \geq -K(0).
\end{equation*}
Taking a limit as $n \rightarrow \infty$ gives:
\begin{equation*}
  -\min\prts{ \frac{1}{\sqrt{V}}, \frac{L_\gen'}{2} } \expect{T_\gen} \geq -K(0),
\end{equation*}
where $\expect{T_\gen} = \sum_{t=0}^\infty \prob{ T_\gen > t }$.  This proves the lemma.  Note that this argument is similar to a drift lemma given in Proposition 5.3 in \cite{Soren:prob}.
\end{IEEEproof}
\end{lemma}

Lemma \ref{lem:transient_gen2} implies that the expected transient time is $O(V^{1.5})$ for problem \eqref{problem:sto} satisfying Assumptions \ref{ass:uniqueness} and \ref{ass:nonpolyhedral}.

\subsection{Convergence time in a steady state}
Once Algorithm \ref{alg:dpp} enters the steady state, the following property holds.

\begin{lemma}
  \label{lem:distance_gen}
  When $V$ is sufficiently larger that $B_\gen(V) < SV$, $B_\gen' \leq SV$, and $\sqrt{V} \geq 2/{L_\gen'}$, under Assumptions \ref{ass:uniqueness} and \ref{ass:nonpolyhedral}, for any time $t \geq T_\gen$, the following holds 
\begin{align}
  \expect{ \norm{ Q(t) - V\dvar^\ast } } & \leq U_\gen(V) \label{eq:ex_gen}\\
  \expect{ \norm{ Q(t) - V\dvar^\ast }^2 } & \leq U_\gen'(V) \label{eq:ex2_gen}
\end{align}
where $U_\gen(V), U_\gen'(V), r_\gen(V), \rho_\gen(V),$ and $D_\gen(V)$ are:
\begin{align*}
  U_\gen(V) &\defeq \frac{ \log\prtr{D_\gen(V) + e^{r_\gen(V) B_\gen(V)}} }{ r_\gen(V) } &&(=O(\sqrt{V})), \\
  U_\gen'(V) &\defeq \frac{ 2\prtr{ D_\gen(V) + e^{r_\gen(V) B_\gen(V)} } }{ r_\gen(V)^2 } &&(=O(V)), \\
  r_\gen(V) &\defeq \frac{3}{6C\sqrt{V} + \sqrt{2C}} &&(= O(\frac{1}{V})), \\
  \rho_\gen(V) &\defeq 1 - \frac{3}{12CV + 2\sqrt{2CV}} &&(= O(1)), \\
  D_\gen(V) &\defeq \frac{e^{r_\gen(V) B_\gen(V)}\prtr{ e^{r_\gen(V) \sqrt{2C}} - \rho_\gen(V)}}{ 1 - \rho_\gen(V) } &&(=O(e^\frac{1}{\sqrt{V}})).
\end{align*}

\begin{IEEEproof}
Recall that $K(t) \defeq \norm{ Q(t) - V\dvar^\ast }$.  From Lemmas \ref{lem:k} and \ref{lem:gen_onestep}, constants in Theorem \ref{thm:concentration} are $\delta = \sqrt{2C}, K = B_\gen(V)$, and $\beta = 1/\sqrt{V}$.  By the definition of $T_\gen$, $k_0 = K(T_\gen) < K$, and Theorem \ref{thm:concentration} implies, for any $t \geq T_\gen$, that
\begin{equation}
  \label{lem_distance_gen:1}
  \expect{ e^{rK(t)} } \leq D + \prtr{ e^{rk_0} - D }\rho^\prtr{t - T_\gen} \leq D + e^{rk_0}.
\end{equation}
Jensen's inequality implies that $e^{r\expect{K(t)}} \leq \expect{ e^{rK(t)} }$, and we have $e^{r\expect{K(t)}} \leq D + e^{rk_0}$.  Taking logarithm and dividing by $r$ proves \eqref{eq:ex_gen}.

Chernoff bound \cite{Ross:Stochastic} implies that, for any $m \in \RealSet_+$,
\begin{equation}
  \label{lem_distance_gen:2}
  \prob{ K(t) \geq m } \leq e^{-rm} \expect{ e^{r K(t) } } \leq \prtr{D + e^{rk_0}} e^{-rm},
\end{equation}
where the last inequality uses \eqref{lem_distance_gen:1}.

Since $K(t)^2$ is always non-negative, it can be shown that $\expect{ K(t)^2 } = 2 \int_0^\infty m \prob{ K(t) \geq m } dm$ by the integration by parts.  Using \eqref{lem_distance_gen:2}, we have
\begin{align*}
  \expect{ K(t)^2 } &\leq 2\prtr{D + e^{rk_0}} \int_0^\infty m e^{-rm} dm.
\end{align*}
Performing the integration by parts proves \eqref{eq:ex2_gen}.
\end{IEEEproof}
\end{lemma}

The convergence results in the steady state are as follows.
\begin{theorem}
  \label{thm:nonpoly}
  When $V$ is sufficiently large that $B_\gen(V) < SV, B_\gen' \leq SV$, and $\sqrt{V} \geq 2/{L_\gen'}$, then for any time $t_0 \geq T_\gen$ and any positive integer $T$, the objective cost converges as
\begin{multline}
  \label{eq:obj_nonpoly}
  \expect{ f( \tavg{x}(t_0, T) ) } - f^\optvar \leq \frac{2M_f U_\gen(V)}{T} \\+ \frac{U_\gen'(V) + 4VU_\gen(V)\norm{\dvar^\ast}}{2TV} + \frac{C}{V}
\end{multline}
and the constraint violation is upper bounded by
\begin{equation}
  \label{eq:con_nonpoly}
  \expect{ g_j(\tavg{x}(t_0, T)) } \leq \frac{2U_\gen(V)}{T} + \frac{2M_{gj} U_\gen(V)}{T}.
\end{equation}
\begin{IEEEproof}
From Theorem \ref{thm:gen_convergence}, the objective cost converges as \eqref{eq:obj_gen}.  Since $T_\poly \leq t_0 < t_0+T$, we use results in Lemma \ref{lem:distance_gen} to upper bound $\expect{ \norm{ Q(t) - V\dvar^\ast } }$ and $\expect{ \norm{ Q(t) - V\dvar^\ast }^2 }$ for $t_0$ and $t_0 + T$.
%\begin{multline}
%  \label{thm-gen:1}
%  \expect{ f( \tavg{x}(t_0, T) ) } - f^\optvar \leq \frac{M_f}{T} \expect{ \norm{ Z(t_0 + T) - Z(t_0) } } \\ + \frac{1}{2TV} \expect{ \norm{Q(t_0)}^2 - \norm{Q(t_0 + T)}^2 } + \frac{C}{V}.
%\end{multline}
Terms in the right-hand-side of \eqref{eq:obj_gen} are bounded by
\begin{align}
  & \expect{ \norm{ Z(t_0 + T) - Z(t_0) } } \leq \expect{ \norm{ Q(t_0 + T) - Q(t_0) } }  \notag\\
  & \quad \leq \expect{ K(t_0 + T) + K(t_0) } \leq 2U_\gen(V). \label{thm-gen:2}
\end{align}
Lemma \ref{lem:Qsquare} implies that
\begin{align}
  &\expect{ \norm{Q(t_0)}^2 - \norm{Q(t_0 + T)}^2 } \notag\\
  & \quad \leq \expect{ K(t_0)^2 + 2\norm{ V\dvar^\ast }\prts{ K(t_0) + K(t_0+T) } } \notag\\
  & \quad \leq U_\gen'(V) + 4VU_\gen(V) \norm{\dvar^\ast}. \label{thm-gen:3}
\end{align}
Substituting bounds \eqref{thm-gen:2} and \eqref{thm-gen:3} into \eqref{eq:obj_gen} proves \eqref{eq:obj_nonpoly}.

The constraint violation converges as \eqref{eq:con_gen} where $T_\poly \leq t_0 < t_0 + T$, and Lemma \ref{lem:distance_gen} can be utilized.
%\begin{multline}
%  \label{thm-gen:4}
%  \expect{ g_j(\tavg{x}(t_0, T)) } \leq \frac{1}{T} \expect{ W_j(t_0 + T) - W_j(t_0) } \\+ \frac{M_{gj}}{T} \expect{ \norm{ Z(t_0+T) - Z(t_0) } }.
%\end{multline}
The last term in the right-hand-side of \eqref{eq:con_gen} is bounded in \eqref{thm-gen:2}.  The first term is bounded by
\begin{align*}
  & \expect{ W_j(t_0 + T ) - W_j(t_0) } \\
  & \quad \leq \expect{ \abs{ W_j(t_0 + T ) - V\dvar^\ast } + \abs{ W_j(t_0) - V\dvar^\ast } } \\ 
  & \quad \leq \expect{ K(t_0 + T) + K(t_0) } \leq 2U_\gen(V).
\end{align*}
Substituting the above bound and \eqref{thm-gen:2} into \eqref{eq:con_gen} proves \eqref{eq:con_nonpoly}.
\end{IEEEproof}
\end{theorem}

The implication of Theorem \ref{thm:nonpoly} is as follows.  When the average starts in the steady state, the deviation from the optimal cost is $O(\sqrt{V}/T + 1/V)$, and the constraint violation is bounded by $O(\sqrt{V}/T)$.  Note that this can be shown by substituting $B_\gen(V), \tilde{r}(V), \tilde{\rho}(V), \tilde{D}_1(V)$, and $\tilde{D}_2(V)$ into \eqref{eq:obj_nonpoly} and \eqref{eq:con_nonpoly}.  By setting $V = 1/\epsilon$, the convergence time is $O(1/\epsilon^{1.5})$.

\section{Simulation}
\label{sec:simulation}

\subsection{Staggered Time Averages}
In order to take advantage of the improved convergence times, computing time averages must be started in the steady state phase shortly after the transient time.  To achieve this performance without knowing the end of the transient phase, time averages can be restarted over successive frames whose frame lengths increase geometrically.  For example, if one triggers a restart at times $2^k$ for integers $k$, then a restart is guaranteed to occur within a factor of 2 of the time of the actual end of the transient phase.

\subsection{Results}
This section illustrates the convergence times of the drift-plus-penalty Algorithm \ref{alg:dpp} under locally-polyhedron and locally-non-polyhedron assumptions.  Let $\setrandom = \prtc{ 0, 1, 2 }, \set{X}_0 = \prtc{ (0,0) }, \set{X}_1 = \prtc{ (-5, 0), (0, 10) }, \set{X}_2 = \prtc{ (0,-10), (5,0) }$, and $(\ssp_0, \ssp_1, \ssp_2) = (0.1, 0.6, 0.3)$.  A formulation is
\begin{align}
  \label{eq:sim_formulation}
  \minimize \quad & \expect{ f(\bar{x}) } \\
  \subjectto \quad & \expect{ 2 \bar{x}_1 + \bar{x}_2 } \geq 1.5 \notag\\
  & \expect{ \bar{x}_1 + 2 \bar{x}_2} \geq 1.5 \notag\\
  & (x_1(t), x_2(t)) \in \set{X}_{\random(t)}, \quad t \in \{0, 1, 2, \dotsc\} \notag
\end{align}
where function $f$ will be given for different cases.

\begin{figure}
  \centering
  \includegraphics[scale=0.40]{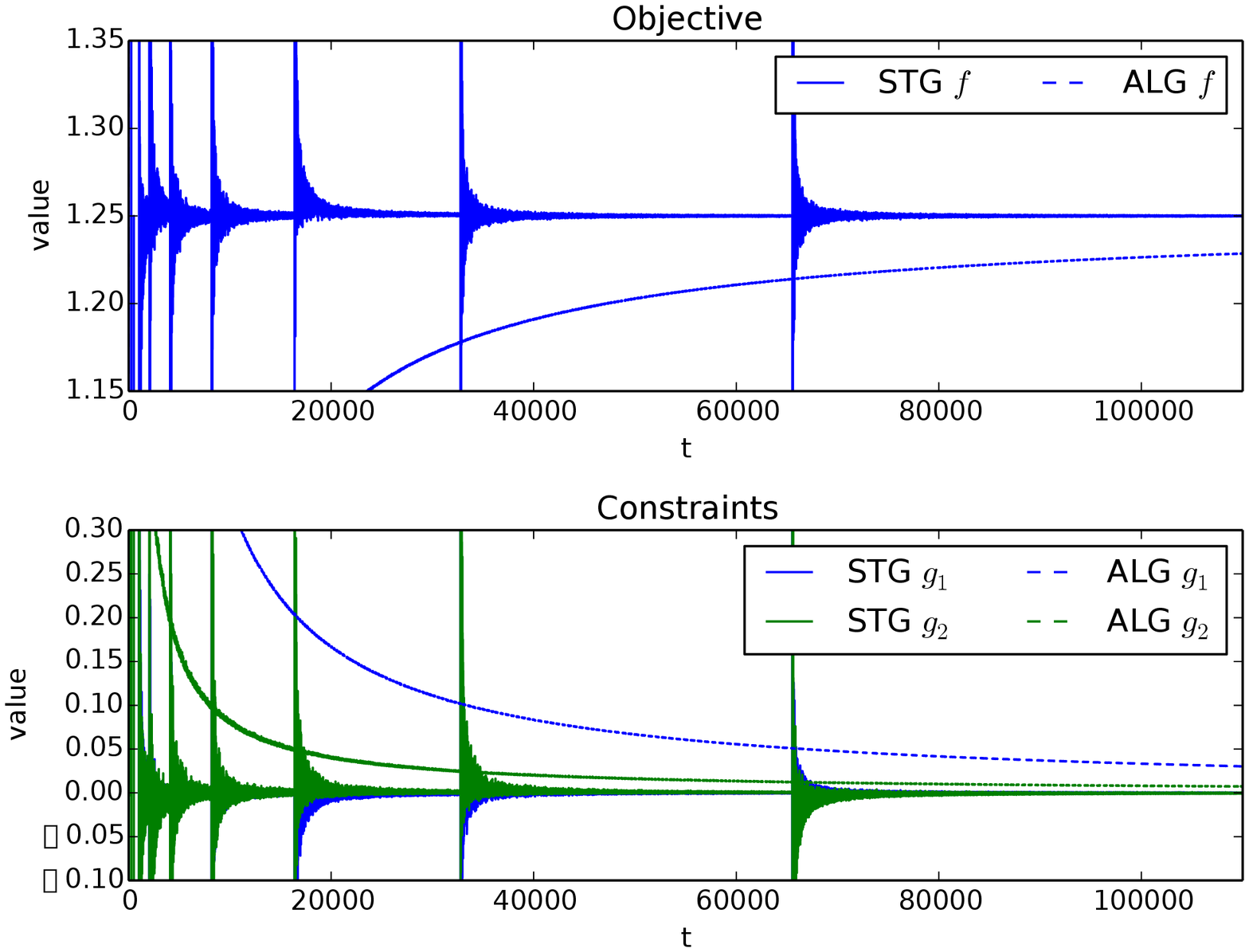}
  \caption{Processes of solving \eqref{eq:sim_formulation} with $f(x) = 1.5 x_1 + x_2$}
  \label{fig:poly_unique}
\end{figure}

Under locally-polyhedron assumption, let $f(x) = 1.5 x_1 + x_2$ be the objective function of problem \eqref{eq:sim_formulation}.  In this setting, the optimal value is $1.25$ where $\tavg{x}_1 = \tavg{x}_2 = 0.5$.  Figure \ref{fig:poly_unique} shows the values of objective and constraint functions of time-averaged solutions.  It is easy to see the improved convergence time $O(1/\epsilon)$ from the staggered time averages (STG) compared to the convergence time $O(1/\epsilon^2)$ of Algorithm \ref{alg:dpp} (ALG).

\begin{figure}
  \centering
  \includegraphics[scale=0.40]{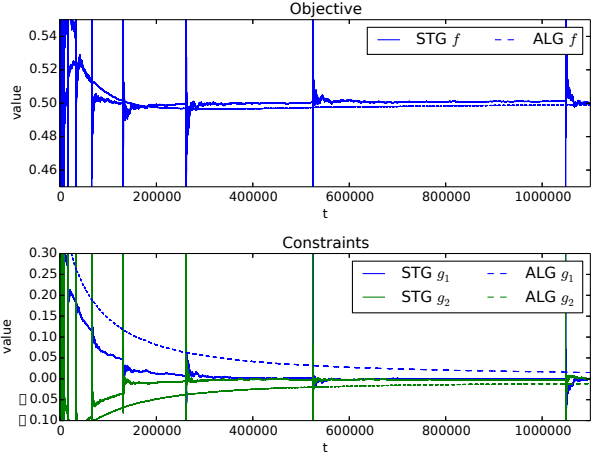}
  \caption{Processes of solving \eqref{eq:sim_formulation} with $f(x) = x_1^2 + x_2^2$}
  \label{fig:nonpoly_unique}
\end{figure}

Under locally-non-polyhedral assumption, let $f(x) = x_1^2 + x_2^2$ be the objective function of problem \eqref{eq:sim_formulation}.  Note that the optimal value of this problem is $0.5$ where $\bar{x}_1 = \bar{x}_2 = 0.5$.  Figure \ref{fig:nonpoly_unique} shows the values of objective and constraint functions of time-averaged solutions.  It can be seen from the constraints plot that the staggered time averages converges faster than Algorithm \ref{alg:dpp}.  This illustrates the different between convergence times $O(1/\epsilon^{1.5})$ and $O(1/{\epsilon^2})$.

\begin{figure}
  \centering
  \includegraphics[scale=0.40]{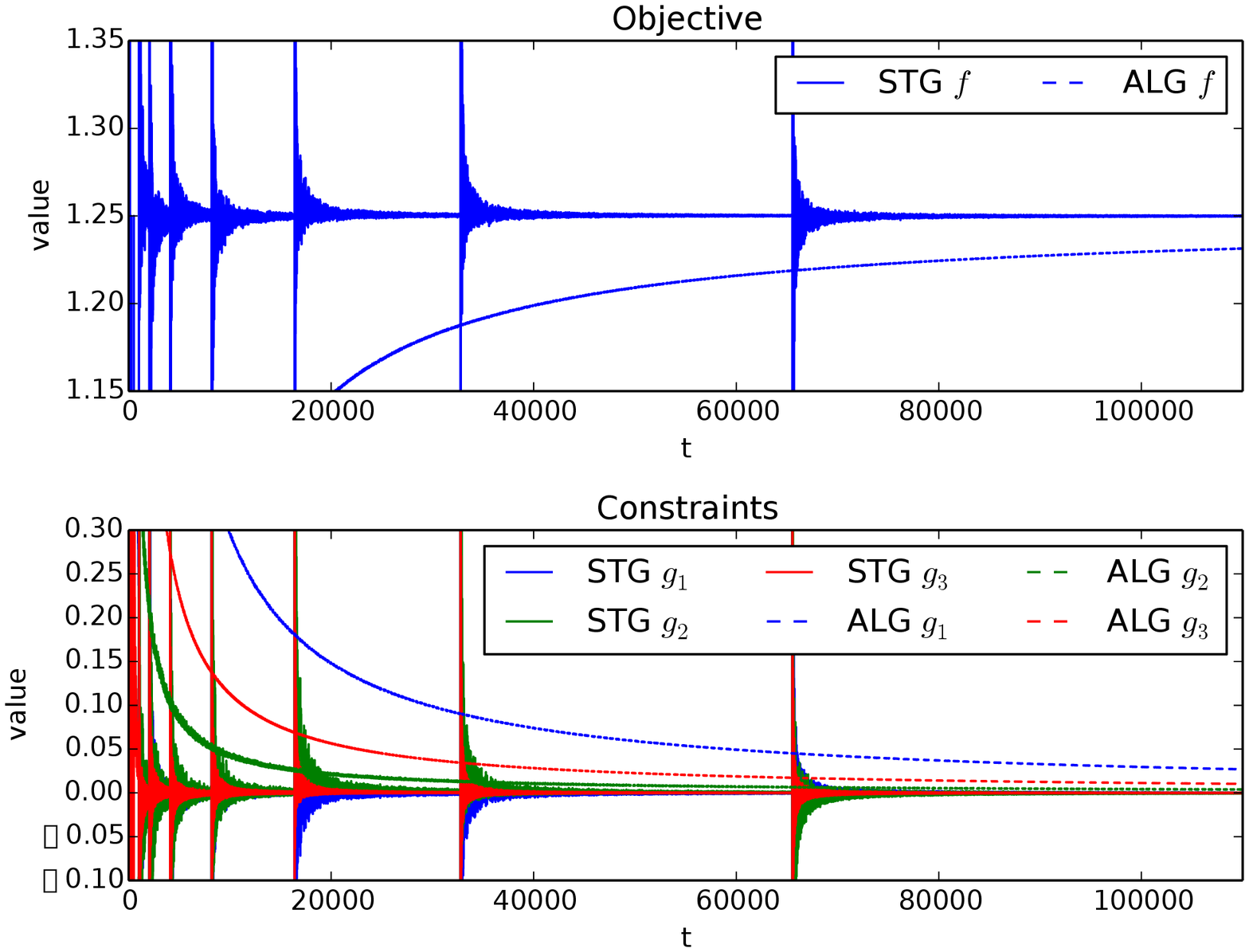}
  \caption{Iterations in processes of solving \eqref{eq:sim_formulation} with $f(x) = 1.5 x_1 + x_2$ and an additional constraint $\expect{\bar{x}_1 + \bar{x}_2} \geq 1$}
  \label{fig:poly_nonunique}
\end{figure}

Figure \ref{fig:poly_nonunique} illustrates the convergence time of problem \ref{eq:sim_formulation} with $f(x) = 1.5 x_1 + x_2$ and additional constraint $\expect{\bar{x}_1 + \bar{x}_2 } \geq 1$.  The dual function of this formulation has non-unique vector of Lagrange multipliers.  The Comparison of Figures \ref{fig:poly_nonunique} and \ref{fig:poly_unique} shows that there is no difference in the order of convergence time.

\begin{figure}
  \centering
  \includegraphics[scale=0.40]{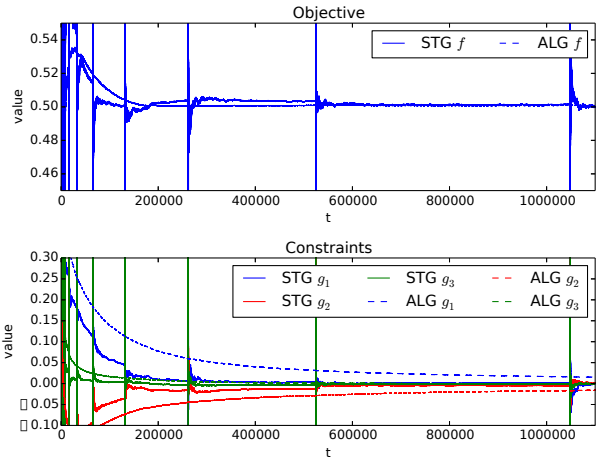}
  \caption{Iterations in processes of solving \eqref{eq:sim_formulation} with $f(x) = x_1^2 + x_2^2$ and an additional constraint $\expect{\bar{x}_1 + \bar{x}_2} \geq 1$}
  \label{fig:nonpoly_nonunique}
\end{figure}

Figure \ref{fig:nonpoly_nonunique} illustrates the convergence time of problem \ref{eq:sim_formulation} with $f(x) = 1.5 x_1 + x_2$ and additional constraint $\expect{ \bar{x}_1 + \bar{x}_2 }\geq 1$.  The dual function of this formulation has non-unique vector of Lagrange multipliers.  The Comparison of Figures \ref{fig:nonpoly_nonunique} and \ref{fig:nonpoly_unique} shows that there is no difference in the order of convergence time.

\section{Conclusion}
\label{sec:conclusion}
We consider the time-average stochastic optimization problem with a non-convex decision set.  The problem can be solved using the drift-plus-penalty algorithm, which converges in $O(1/\epsilon^2)$.  After we analyze the transient and steady state phases of the algorithm, the convergence time can be improved by performing time average in the steady state.  we prove that the improved convergence time is $O(1/\epsilon)$ under the locally-polyhedral assumption and is $O(1/\epsilon^{1.5})$ under the locally-non-polyhedral assumption.

\section*{Appendix}
The convergence time of the objective cost as a function of vectors $y(t)$ is firstly proven.  This proof requires the following theorem proven in \cite{Neely:SNO}.

\begin{theorem}
  \label{thm:random_policy}
  There exists a randomized policy $x^\diamond(t), y^\diamond(t)$ that only depends on $\random(t)$ such that for all $t \in \prtc{0, 1, 2, \dotsc}$:
\begin{align*}
  \expect{ f(y^\diamond(t) } & = f^\optvar \\
  \expect{ g(y^\diamond(t)) } & \leq 0 \\
  \expect{ x^\diamond(t)} & = \expect{ y^\diamond(t) } \\
  \prtr{ x^\diamond(t), y^\diamond(t) } & \in \set{X}_{\random(t)} \times \set{Y}.
\end{align*}
\end{theorem}

\begin{lemma}
  \label{lem:objective_y}
  Let $\prtc{ (x(t), y(t), W(t), Z(t)) }_{t=0}^\infty$ be a sequence generated by Algorithm \ref{alg:dpp}.  For any positive integer $T$ and starting time $t_0$, it holds that
\begin{multline*}
  \expect{ f( \tavg{y}(t_0, T) } - f^\optvar \leq \\\frac{C}{V} + \frac{1}{2TV} \expect{ \norm{Q(t_0)}^2 - \norm{Q(t_0 + T)}^2 }.
\end{multline*}

\begin{IEEEproof}
Since decision $(x(t), y(t))$ minimizes the right-hand-side of \eqref{eq:dpp_bound}, the following holds for any other decisions including the randomized policy in Theorem \ref{thm:random_policy}:
\begin{align*}
  \drift(t) + V f(y(t)) & \leq C + Q(t)^\tr h(x(t), y(t)) + V f(y(t)) \\
  & \leq C + Q(t)^\tr h(x^\diamond(t), y^\diamond(t) + V f(y^\diamond(t)).
\end{align*}
Taking the conditional expectation of the above bound gives
\begin{align*}
  & \expect{ \drift(t) | Q(t) } + V \expect{ f(y(t)) | Q(t) } \\
  & \quad \leq C + Q(t)^\tr \expect{ h(x^\diamond(t), y^\diamond(t) | Q(t) } + V \expect{ f(y^\diamond(t)) | Q(t) } \\
  & \quad \leq C + V f^\optvar,
\end{align*}
where the last inequality uses properties of the randomized policy in Theorem \ref{thm:random_policy}.

Taking expectation and using iterated expectation leads to
\begin{equation*}
  \frac{1}{2} \expect{ \norm{Q(t+1)}^2 - \norm{Q(t)}^2 } + V \expect{ f(y(t)) } \leq C + V f^\optvar.
\end{equation*}
Summing from $t = t_0$ to $t_0 + T-1$ yields
\begin{multline*}
  \frac{1}{2} \expect{ \norm{Q(t_0 + T)}^2 - \norm{Q(t_0)}^2 } + V \expect{ \sum_{t=t_0}^{t_0+T-1} f(y(t)) } \\\leq CT + VT f^\optvar.
\end{multline*}
Dividing by $T$, using the convexity of $f$, and rearranging terms proves the lemma.
\end{IEEEproof}
\end{lemma}

The constraint violation as a function of vectors $y(t)$ is the following.
\begin{lemma}
  \label{lem:constraint_y}
  Let $\prtc{ (x(t), y(t), W(t), Z(t)) }_{t=0}^\infty$ be a sequence generated by Algorithm \ref{alg:dpp}.  For any positive integer $T$ and starting time $t_0$, it holds for $j \in \prtc{1, \dotsc, J}$ that
\begin{equation*}
  \expect{ g_j( \tavg{y}(t_0, T) ) } \leq \frac{1}{T} \expect{ W_j(t_0 + T) - W_j(t_0) }.
\end{equation*}

\begin{IEEEproof}
Dynamic \eqref{eq:W_update} implies that $W_j(t+1) \geq W_j(t) + g_j(y(t))$.  Taking expectation gives $\expect{W_j(t+1)} \geq \expect{W_j(t)} + \expect{g_j(y(t))}$.  Summing from $t=t_0$ to $t_0 + T-1$ yields $\expect{W_j(t_0 + T)} \geq \expect{W_j(t_0)} + \expect{\sum_{t=t_0}^{t_0 + T-1} g_j(y(t)) }$. Dividing by $T$, using the convexity of $g_j$, and rearranging terms proves the lemma.
\end{IEEEproof}
\end{lemma}

The following result is used to translates the results in Lemmas \ref{lem:objective_y} and \ref{lem:constraint_y} to the bounds as functions of vectors $x(t)$.
\begin{lemma}
  \label{lem:y_to_x}
  Let $\prtc{ (x(t), y(t), W(t), Z(t)) }_{t=0}^\infty$ be a sequence generated by Algorithm \ref{alg:dpp}.  For any positive integer $T$ and starting time $t_0$, it holds that
\begin{equation*}
  \tavg{x}(t_0, T) - \tavg{y}(t_0, T) = \frac{1}{T} \prts{ Z(t_0 + T) - Z(t_0) }.
\end{equation*}

\begin{IEEEproof}
Dynamic \eqref{eq:Z_update} implies that $Z(t+1) = Z(t) + \prts{x(t) - y(t)}$.  Summing from $t=t_0$ to $t_0 + T-1$ yields $Z(t_0 + T) = Z(t_0) + \sum_{t=t_0}^{t_0 + T-1} \prts{x(t) - y(t)}$. Dividing by $T$ and rearranging terms proves the lemma.
\end{IEEEproof}
\end{lemma}

Finally, Lemma \ref{thm:gen_convergence} is proven.
%\begin{theorem}
%  \label{thm:gen_convergence}
%  Let $\prtc{ (x(t), y(t), W(t), Z(t)) }_{t=0}^\infty$ be a sequence generated by Algorithm \ref{alg:dpp}.  For any positive integer $T$ and starting time $t_0$, the objective cost converges as
%\begin{multline}
%  \label{eq:obj_gen}
%  \expect{ f( \tavg{x}(T) ) } - f^\optvar \leq \frac{C}{V} + \frac{1}{2TV} \expect{ \norm{Q(0)}^2 - \norm{Q(T)}^2 } \\+ \frac{M_f}{T} \expect{ \norm{ Z(T) - Z(0) } },
%\end{multline}
%and the constraint violation for every $j \in \prtc{1, \dotsc, J}$ is
%\begin{multline}
%  \label{eq:con_gen}
%  \expect{ g_j(\tavg{x}(T)) } \leq \frac{1}{T} \expect{ W_j(T) - W_j(0) } \\+ \frac{M_{gj}}{T} \expect{ \norm{ Z(T) - Z(0) } }.
%\end{multline}

\begin{IEEEproof}
For the first part, Lipschitz continuity \eqref{eq:f_lip} implies
\begin{multline*}
  f( \tavg{x}(t_0, T) ) - f^\optvar \leq f( \tavg{y}(t_0, T) ) - f^\optvar \\+ M_f \norm{ \tavg{x}(t_0, T) - \tavg{y}(t_0, T) }.
\end{multline*}
Taking expectation yields
\begin{multline*}
  \expect{ f( \tavg{x}(t_0, T) ) } - f^\optvar \leq \expect{ f( \tavg{y}(t_0, T) ) } - f^\optvar \\+ M_f \expect{ \norm{ \tavg{x}(t_0, T) - \tavg{y}(t_0, T) } }.  
\end{multline*}
Applying the results in Lemmas \ref{lem:objective_y} and \ref{lem:y_to_x} proves \eqref{eq:obj_gen}.

For the last part, Lipschitz continuity \eqref{eq:gj_lip} implies that
\begin{multline*}
  g_j( \tavg{x}(t_0, T) ) \leq g_j( \tavg{y}(t_0, T) ) + M_{gj} \norm{ \tavg{x}(t_0, T) - \tavg{y}(t_0, T) } \\ j \in \prtc{1, \dotsc, J}.
\end{multline*}
Taking expectation yields
\begin{multline*}
  \expect{ g_j( \tavg{x}(t_0, T) ) } \leq \expect{ g_j( \tavg{y}(t_0, T) ) } \\+ M_{gj} \expect{ \norm{ \tavg{x}(t_0, T) - \tavg{y}(t_0, T) } }.
\end{multline*}
Applying the results in Lemmas \ref{lem:constraint_y} and \ref{lem:y_to_x} proves \eqref{eq:con_gen}.
\end{IEEEproof}

\bibliographystyle{IEEEtran}
\bibliography{Reference}

\end{document}